\documentclass[12pt]{article}
\usepackage[applemac]{inputenc}
\usepackage{amsmath,amssymb}
\usepackage[T1]{fontenc}
\usepackage{amsmath}
\usepackage{amssymb}
\usepackage{amstext}
\input amssym.def
\input amssym.tex
\usepackage{epsfig}
\usepackage[]{graphicx}
\usepackage{color}
\usepackage{tikz}%
\providecommand{\U}[1]{\protect\rule{.1in}{.1in}}
\usetikzlibrary{arrows,shapes,automata,petri}

\newcommand{\dproof}{\noindent {Proof.} \quad}
\newcommand{\fproof}{\hfill $\square$ \bigskip}
\newtheorem{definition}{Definition}[section]
\newtheorem{example}{Example}[section]
\newtheorem{theorem}[definition]{Theorem}

\newtheorem{remark}[definition]{ \it Remark}

\newtheorem{lemma}[definition]{Lemma}

\numberwithin{equation}{section}

\def\RR{{\mathbb{ R}}}

\def\RR{{\mathbb{R}}}

\def\1B{\text{1\!\!I}}

\def\P{\mathbb{P}}
\def\E{\mathbb{E}}
\def\R{\mathbb{R}}

\author{\bf Jasmina \DJ or\dj evi\'{c} \thanks{  {\it Address:} Faculty of Science and Mathematics, University of Ni\v s, Vi\v segradska 33, 18000 Ni\v s, Serbia,      \newline    {\it E-mail:}   jasmina.djordjevic@pmf.edu.rs,, djordjevichristina@gmail.com\newline Supported by the Ministry of Science, Technological Development and Innovation of Republic of Serbia, agreement no.  451-03-34/2026-03/200124} \bf Bernt \O ksendal
\thanks{ {\it Address:} Department of Mathematics, University of Oslo, Blindern 0316 Oslo, Norway,       \newline    {\it E-mail:}  oksendal@math.uio.no.}
}

\date{7 April 2026}
\title{Conditional stochastic differential equations driven by fractional Brownian motion}

\begin{document}
\maketitle

\begin{abstract}
The aim of this paper is to analyse a WIS-stochastic differential equation driven by fractional Brownian motion with $H>\tfrac{1}{2}$. For this, we  summarise the theory of fractional white noise and prove a fundamental  $L^2$-estimate for WIS-integrals. We apply this to prove the existence and uniqueness of a solution in $L^2(P)$ of a conditional WIS-stochastic differential equation driven by a fractional Brownian motion with $H>\tfrac{1}{2}$ under Lipschitz  conditions on its coefficients.

\bigskip


\end{abstract}

\section{Introduction}
Here is a motivation for the study in this paper:\\

Suppose we initially model the growth of a (biological, for example) system in a noisy environment by a classical It\^o type stochastic differential equation (SDE)  of the form 
\begin{align}
dX_t= \alpha (t,X_t) dt + \beta (t,X_t) dB(t); t \geq 0
\end{align}
where $B(t)=B(t,\omega)$ is a Brownian motion, and for all $x$ the processes 
$t \mapsto \alpha(t,x,\omega)$ and
$t \mapsto  \beta(t,x,\omega)$
are both adapted to the filtration $\mathcal{F}$  of $B(\cdot)$.
Using white noise notation  this can be written
\begin{align}
\tfrac{dX_t}{dt}= \alpha(t,X_t)+ \beta(t,X_t) \diamond \overset{\bullet}{B}(t),
\end{align}
where $\diamond$ denotes Wick product.

Now suppose that in addition to the rapidly fluctuating  noise $ \overset{\bullet}{B}(t)$ there is a slower fluctuating, memory noise $\overset{\bullet}{B}^{(H)}(t)$ in the system, where $B^{(H)}(t)$ is a \emph{fractional} Brownian motion defined on the same probability space as $B(\cdot)$, and with Hurst coefficient $H > \tfrac{1}{2}$.
Then a model for the dynamics of the system could be a \emph{conditional} stochastic differential equation of the form
\begin{align} \label{sde1}
dX_t &= \alpha(t,X_t) dt +\beta(t,\E[X_t | \mathcal{F}_t]) dB(t) + \sigma(t,\E[X_t | \mathcal{F}_t] )dB^{(H)}(t), 
\end{align}
or, equivalently,
\begin{align} \label{sde1+}
dX_t &= \alpha(t,X_t) dt +\beta(t,\E[X_t | \mathcal{F}_t]) dB(t) + \sigma(t,\E[X_t | \mathcal{F}_t] )dB^{(H)}(t), \\
&\text{ where}\nonumber\\
Y_t&=\E[X_t | \mathcal{F}_t]
\end{align}
is the conditional expectation of $X_t$ given the history $\mathcal{F}_t$, defined to be the $\sigma$-algebra generated by $B(s,\omega); s \leq t$. See Section 2 for more details about the construction of $B^{(H)}(t)$ and its relation to $B(t)$.

In this paper we study existence and uniqueness of solutions $X_t=X(t)=X(t,\omega);t\geq 0,\omega \in \Omega$  of such conditional stochastic differential equations,  where $X_0=Z$ is a given random variable in $L^2(\P)$, independent of both $B(\cdot)$ and the fractional Brownian motion (fBm for short)  $B^{(H)}(\cdot)$. We assume that the Hurst coefficient $H$ of $B^{(H)}$ is in the interval $(\tfrac{1}{2},1)$ and the stochastic integrals are interpreted as Wick-Ito-Skorohod (WIS)-integrals, i.e. the SDE \eqref{sde1} (suppressing the $\omega$) can be written as the following equation in the Hida stochastic distribution space $(\mathcal{S})^{*}$:
\begin{align}
\tfrac{d}{dt}X_t&=b(t,X_t) +\beta(t,\E[X_t | \mathcal{F}_t]) \diamond \overset{\bullet}{B}(t)+ \sigma (t,\E[X_t | \mathcal{F}_t]) \diamond \overset{\bullet}{B}^{(H)}_t;\quad 0 \leq t \leq T\\
X(0)&=Z,\nonumber
\end{align}
where $\diamond$ denotes the Wick product, and $\overset{\bullet}{B}_t^{(H)}=\tfrac{d}{dt} B_t^{(H)}$ is the fractional white noise. 
We do not assume that the coefficients are adapted to the filtration generated by $B^{(H)}(\cdot)$.

A general existence of uniqueness result for (unconditional)  WIS-type SDEs driven by fractional Brownian motion is obtained by Yaozhong Hu \cite{Hu}. He proves the existence and uniqueness of a solution up to a certain stochastic time, under (quite restrictive) conditions  on the coefficients, involving their Malliavin derivatives . 

We prove that under condition of Lipschitz continuity of the coefficients there is a unique measurable solution $X(t)$ of \eqref{sde1} with $\E[\int_0^T X^2(t)dt] < \infty$ (Theorem \ref{sde3}).
The proof is based on an $L^2$-estimate of WIS-integrals of independent interest (Theorem 3.1). The above stochastic system combines several advanced features arising from the interaction of mixed noise, partial information, and mean-field type dependence.

First, the state dynamics are driven by both a standard Brownian motion and a fractional Brownian motion with Hurst parameter $H>1/2$. While the Brownian motion represents a classical source of randomness with independent increments, the fractional Brownian motion introduces long-range dependence and persistent behavior. In particular, the presence of the fractional component destroys the Markov property and prevents the use of classical It\^o calculus in a straightforward way.

Second, the filtration $\{\mathcal{F}_t\}_{t\geq 0}$ is generated solely by the Brownian motion $B(s); s \leq t$. Consequently, the fractional Brownian motion is not adapted to this filtration and can be interpreted as an unobservable or latent noise source. From the point of view of the available information, the system is therefore only partially observed.

The process
\[
Y_t = \mathbb{E}[X_t \mid \mathcal{F}_t]
\]
naturally appears as the filtered estimate (or projection) of the state process $X_t$ onto the observable filtration. In this sense, $Y_t$ represents the best mean-square estimate of $X_t$ given the available information. This introduces an intrinsic connection with nonlinear filtering theory, where one studies the evolution of conditional distributions or conditional expectations under incomplete information.

Moreover, the coefficients $\beta$ and $\sigma$ depend on the filtered process $Y_t$ rather than directly on $X_t$. This creates a nontrivial feedback loop: the state process $X_t$ depends on $Y_t$, while $Y_t$ itself is defined in terms of $X_t$ via conditional expectation. As a result, the system can be interpreted as a conditional mean-field (or McKean--Vlasov type) equation under partial information, where the interaction is mediated through the observable component of the system.    Conditional expectation: Since $B^H$ is independent of $\mathcal{F}_t$, its contribution averages out in $Y_t$:
    \[
    Y_t = \mathbb{E}[X_t \mid \mathcal{F}_t],
    \]
representing the best mean-square estimate of $X_t$ given the available information. From an applied perspective, such models arise naturally in situations where part of the randomness is not directly observable. The fractional Brownian motion may represent hidden environmental effects, memory, or long-term correlations (for instance in biological, epidemiological, or financial systems), while the Brownian motion represents short-term observable fluctuations. The filtered process $Y_t$ then corresponds to the real-time estimate of the underlying state based on observable data.

Finally, the existence and uniqueness result can be interpreted as a well-posedness statement for a stochastic system with hidden memory effects and feedback through filtering.

 As far as we know this estimate and the resulting theorem are new. 

\subsection{Some historical remarks and earlier literature}
Mandelbrot and Van Ness \cite{MN} introduced the concept of fractional Brownian motion as a generalisation of classical Brownian motion. They defined fBm as a Gaussian process with stationary increments and a specific covariance structure that depends on a parameter called the Hurst exponent. They demonstrated that if $H>\tfrac{1}{2}$ then  fBm exhibits long-range dependence, which means that its increments are correlated over a wide range of time scales.  With decrease of Hurst parameter, trajectory of fBm is more unstable (values are changing quicker) and vice versa. This property and some kind of control with respect to Hurst parameter provide wide application of fBm in real problems. There are also interesting connections to models of anomalous diffusion. See e.g. Eliazar \& Kachman \cite {EK} and the references therein.

Fractional Brownian motion is a continuous-time stochastic process with stationary increments and self-similarity, characterised by the Hurst exponent, denoted by $H$.
For $H>\tfrac{1}{2}$ the stochastic integral of a stochastic process $\gamma(s,\omega)$ with respect to $ B^{(H)}(s)$ , denoted by 
$\int_0^t \gamma(s,\omega) dB^{(H)}_s$,  can be defined in (at least) two different ways:
\begin{itemize}
\item
As a pathwise (or forward) integral, 
or 
\item
 as a Wick-Ito-Skorohod (WIS)  integral.
 \end{itemize}
In this paper we concentrate on the second interpretation, which in terms of white noise calculus can be expressed as follows:
\begin{align}
\int_0^t \gamma(s,\omega) dB^{(H)}_s = \int_0^t \gamma(s,\omega) \diamond \tfrac{d}{ds}B^{(H)}_s ds= \int_0^t \gamma(s,\omega)\overset{\bullet}{B}^{(H)}_s ds.
\end{align}

In this sense WIS-type stochastic differential equations (SDEs) driven by fractional Brownian motion (fBm) can be regarded as a natural generalisation of classical It\^o type SDEs driven by classical Brownian motion (corresponding to $H=\tfrac{1}{2}$).  For more information about fBm and its calculus we refer to \cite{BHOZ}, \cite{M} and \cite{JM}.\\

The study of SDEs with fBm finds applications in various fields, including finance, physics, biology, and engineering, where modelling systems with long-range dependence and memory effects is crucial.  Even though the theory of SDEs driven by fBm is still an active area of research, and closed-form solutions are often difficult to obtain, they have been and still are widely explored. In the sequel we briefly recall what is done in this area.

\medskip

The paper \cite{NR} deals with SDEs driven by fBm like our paper, but with a different type of stochastic integration. Therefore their resuls are not comparable to ours. 
The results of \cite{GN} differ from ours as the authors analyse the  stochastic differential equation  driven simultaneously by a multidimensional fractional Brownian motion with Hurst parameter $H > \tfrac{1}{2}$  with a pathwise (forward) interpretation of the integral, and a multidimensional standard Brownian motion. The idea of paper \cite{TV} differs from our case as the authors analyse  linear stochastic evolution equations driven by infinite-dimensional fBm. In \cite{BCO} the authors prove the existence and continuous dependence of mild solutions to stochastic differential equations with non-instantaneous impulses driven by fBm.  In  \cite{JM} using fractional calculus authors prove anticipating Girsanov transformation for fBm and used this to prove well-posedness for SDE with fBm.

Other  problems regarding fBm are analyzed in  \cite{CC},  \cite{DU}, \cite{NS}. 
In \cite{Bu}  a general study of SDEs with non -adapted coefficients are given, using a Skorogod interpretation of the stochastic integral. This could in principle be relevant for our study, if the fractional Brownian motion had been defined via a Volterra type representation of the Brownian motion.  But this is not the definition of fBm used in our paper.


Recently, authors have analysed different problems with fBm. For example, in \cite{CZQ} authors study stochastic harvesting population problems for systems with fBm and formulate control problem for it. In \cite{JME} authors consider fuzzy stochastic differential equations driven by fBm. Further,  in \cite{SA} the author proves the  existence  and  uniqueness  of  mild  solutions  to  stochastic evolution equations with time delay driven by fBm, but with deterministic diffusion coefficients. We also mention that properties of McKean-Vlasov SDEs of hybrid type (driven by both Bm and fBm, as (\ref{sde1}) ) have recently been studied by T. Zhang et al. See [31] and the references therein.

To the best of our knowledge none of the above works cover the results of our paper.

\medskip

Our paper is organised in the following way: 
\begin{itemize}
\item
In  Section 2 we introduce and summarise the theory of fractional white noise and its stochastic calculus. 
\item
 In Section 3  we prove a fundamental  $L^2(\P)-$estimate for WIS-integrals. 
 \item
   Finally, in Section 4 we prove our main result, namely the existence and uniqueness of a solution in $L^2(\P)$ of a WIS-stochastic differential equation driven by a fractional Brownian motion with $H>\tfrac{1}{2}$ under Lipschitz  conditions on its coefficients.
\end{itemize}

\section{Fractional Brownian motion (fBm) and its calculus}

\subsection{The white noise probability space and Brownian motion}
In this subsection we briefly review the basic notation and results of the white noise probability space and its associated Brownian motion. 
Let $\mathcal{S}(\RR)$ denote the Schwartz space of rapidly decreasing smooth functions on $\mathbb{R}$ and define $\Omega={\cal S}'(\RR)$ to be its dual, which is the Schwartz space of tempered distributions, equipped
with the weak-star topology. This space will be the base of our basic \emph{white noise
probability space}, which we recall in the following:
\vskip 0.3cm
As events we will use the family ${\cal
B}({\cal S}'(\RR))$ of Borel subsets of ${\cal S}'(\RR)$, and our
probability measure $\P$ is defined as follows:

\begin{theorem}{\bf (The Bochner--Minlos theorem)}\\
There exists a unique probability measure $\P$ on ${\cal B}({\cal S}'(\R))$
with the following property:
$$\E[e^{i\langle\cdot,\phi\rangle}]:=\int\limits_{\cal S'}e^{i\langle\omega,
\phi\rangle}d\mu(\omega)=e^{-\tfrac{1}{2} \Vert\phi\Vert^2};\quad i=\sqrt{-1}$$
for all $\phi\in{\cal S}(\R)$, where
$\Vert\phi\Vert^2=\Vert\phi\Vert^2_{L^2(\RR)},\quad\langle\omega,\phi\rangle=
\omega(\phi)$ is the action of $\omega\in{\cal S}'(\R)$ on
$\phi\in{\cal S}(\R)$ and $\E=\E_{\P}$ denotes the expectation
with respect to  $\P$.
\end{theorem} 

We will call the triplet $({\cal S}'(\R),{\cal
B}({\cal S}'(\R)),\P)$ the {\it  white noise probability
space\/}, and $\P$ is called the {\it white noise probability measure}.

The measure $\P$ is also often called the (normalised) {\it Gaussian
measure\/} on ${\cal S}'(\R)$. It is not difficult to prove that if $\phi\in L^2(\R)$ and
we choose $\phi_k\in{\cal S}(\R)$ such that $\phi_k\to\phi$ in $L^2(\R)$,
then
$$\langle\omega,\phi\rangle:=\lim\limits_{k\to\infty}\langle\omega,\phi_k\rangle
\quad\text{exists in}\quad L^2(\P)$$
and is independent of the choice of $\{\phi_k\}$. In
particular, if we define
$$\widetilde{B}(t):=\widetilde{B}(t,\omega)=\langle\omega,\chi_
{[0,t]}\rangle; $$ 
then $\widetilde{B}(t,\omega)$ has a $t$-continuous version $B(t,\omega)$, which
becomes a \emph{Brownian motion}, in the following sense:
 
By a \emph{Brownian motion} we mean a family
$\{X(t,\cdot)\}_{t\in\R}$ of random variables on a probability space
$(\Omega,{\cal F},\P)$ such that
\begin{itemize} 
\item
$X(0,\cdot)=0\quad\text{almost surely with respect to } \P,$
\item
$\{X(t,\omega)\}_{t\geq 0,\omega \in \Omega}$ is a continuous and Gaussian stochastic process 
\item
For all $s,t \geq 0$,
$\E[X(s,\cdot) X(t,\cdot)]=\min(s,t).$ 
\end{itemize} 

It can be proved that the process $\widetilde{B}(t,\omega)$ defined above has a modification $B(t,\omega)$ which satisfies all these properties.
This process $B(t,\omega)$ then becomes a Brownian motion. This is the version of Brownian motion we will use to construct the associated \emph{fractional} Brownian motion (fBm), as explained in the next subsection:.

\subsection{The operator $M$}
\begin{definition} \label{fBm}
By a \emph{fractional Brownian motion (fBm)} with Hurst parameter $H \in (0,1)$ we mean a continuous Gaussian stochastic process $\{Y(t)\}_{t\geq 0}$ with mean 0 and covariance
$$\E[Y(s) Y(t)]=C(|s|^{2H} + |t|^{2H}-|t-s|^{2H} ); \text{ for all } s,t\geq 0,$$
for some constant $C>0$.
We will in our paper assume that
\begin{equation}
C=C_H=[2\Gamma(H-\frac{1}{2}) \cos(\frac{\pi}{2}(H-\frac{1}{2}))] ^{-1},
\end{equation}
where $\Gamma$ is the classical Gamma function. See \eqref{MH}.
\
\end{definition}

There are several ways to construct an fBm. We will follow the approach used by, among others,  Elliot and van der Hoek \cite{EvdH}, where fBm is constructed by applying an operator $M$ to the classical Bm. See also \cite{L}, \cite{ST} and \cite{Be1},\cite{Be2}. This approach has the advantage that all the different fBm's corresponding to different Hurst coefficients $H \in (0,1)$ are defined on the same white noise probability space $(\Omega, \mathbb{F},\P)$ introduced below. Moreover, it has several technical advantages, as we have exploited in the proof of our main result (Theorem 3.1)\\

In the following we let $f: \mathbb{R} \mapsto \mathbb{R}$ be a function such that its Fourier transform
\begin{equation}
F(f)(y):= \widehat{f}(y):= \frac{1}{\sqrt{2 \pi}} \int_{\mathbb{R}} e^{-ixy}f(x)dx; \quad y\in \mathbb{R}
 \end{equation}
exists.
Let $F^{-1}$ denote the inverse Fourier transform, defined by
 \begin{align}
 F^{-1} g(x):= \frac{1}{\sqrt{2 \pi}}\int_{\mathbb{R}} e^{ixy} g(y) dy,
 \end{align}
for all integrable functions $g$.

\begin{definition}(\cite{EvdH},(A10))
Let $H \in (0,1)$ and let $\mathcal{S}$ be the Schwartz space of rapidly decreasing smooth functions on $\mathbb{R}$. The operator $M=M^{(H)}: \mathcal{S} \mapsto \mathcal{S}$ is defined by
\begin{align}
Mf(x)&=F^{-1}(|y|^{\frac{1}{2}-H} (Ff)(y))(x)\nonumber\\
&=
\textstyle{\frac{1}{\sqrt{2 \pi}}\int_{\mathbb{R}} e^{ixy} (|y|^{\frac{1}{2} - H}\frac{1}{\sqrt{2 \pi}} \int_{\mathbb{R}} e^{-izy}f(z)dz) dy}
\end{align}
If $H > \frac{1}{2} $ then the operator $M$ can also be defined as follows:
\begin{align}
M_Hf(x)= C_H \int_{\mathbb{R}} \frac{f(y+x)}{| y | ^{\frac{3}{2} - H}} dy; \quad f \in \mathcal{S}.\label{MH}
\end{align}

\end{definition}
\emph{From now on we assume that $H \in (\frac{1}{2},1)$.}\\

Since $H \in (\frac{1}{2},1)$ is fixed, we will suppress the index $H$ and write $M^{(H)} = M$ in the following.\\
We can in a natural way extend the operator $M$  from $\mathcal{S}$ to the space
\begin{align}\label{L2H}
L_{H}^2(\mathbb{R}):= \{ f: \mathbb{R} \mapsto \mathbb{R}; Mf \in L^2(\mathbb{R} )\}.
\end{align}
Note that $L_{H}^2 (\mathbb{R})$ is a Hilbert space when equipped with the inner product
\begin{equation}
(f,g)_H := (f,g)_{L_H^2(\mathbb{R})}=(Mf,Mg)_{L^2 (\mathbb{R})};\quad f,g \in L^2_{H}(\mathbb{R})
\end{equation}
and the corresponding norm
\begin{equation}
\|f\|_H=\Big((Mf,Mf)_{L^2 (\mathbb{R})}\Big)^{\tfrac{1}{2}};\quad f \in L^2_{H}(\mathbb{R}).
\end{equation}

We will need the following results about the operator $M$:
\begin{lemma}
Let $\varphi \in L_H^2(\mathbb{R})$ 
and $\psi \in L^2_{H}(\mathbb{R}), \varphi \geq 0, \psi \geq 0$. Then
\begin{align}
(\varphi,M\psi)_{L^2(\mathbb{R})}=(M\varphi,\psi)_{L^2(\mathbb{R})}
\end{align}
\end{lemma}
\dproof
Since
\begin{equation}
Mf(x)=C_H\int_{\mathbb{R}}|y|^{H-\frac{3}{2}}f(x+y)dy\label{M-add}
\end{equation}
we get by the Fubini-Tonelli theorem
\begin{align*}
(\varphi, M\psi)_{L^2(\mathbb{R})}&= \int_{\mathbb{R}} \varphi(x)M\psi(x)dx= \int_{\mathbb{R}}\varphi(x) C_H \int_{\mathbb{R}}|y|^{H-\frac{3}{2}}\psi(x+y) dy dx\\
&= C_H \int_{\mathbb{R}}\int_{\mathbb{R}}\varphi(z-y) |y|^{H-\frac{3}{2}}\psi(z) dy dz\\
&= C_H \int_{\mathbb{R}}\int_{\mathbb{R}}\varphi(z+y) |y|^{H-\frac{3}{2}}\psi(z) dy dz=(M\varphi,
\psi)_{L^2(\mathbb{R})},
\end{align*}
which completes the proof.
\fproof

\begin{lemma}\label{l1} 
The operator $M$ satisfies following property:
\begin{align} \label{eq3.5}
M^2f(x)= C_H^2 \int_{\mathbb{R}} \int_{\mathbb{R}} |yz|^{H-\frac{3}{2}} f(y+z+x)dy dz; \quad f \geq 0.
\end{align}
\end{lemma}

\dproof 
Since
$$Mf(x)=C_H\int_{\mathbb{R}}|y|^{H-\frac{3}{2}}f(x+y)dy$$
we get
\begin{align*}
M^2f(x)&=M(Mf)(x)= C_H\int_{\mathbb{R}}|z|^{H-\frac{3}{2}} Mf(x+z) dz \\
&= C_H\int_{\mathbb{R}} |z|^{H-\frac{3}{2}} C_H\int_{\mathbb{R}}|y|^{H-\frac{3}{2}}f(x+z+y)dy dz\\
&=C_H^2 \int_{\mathbb{R}} \int_{\mathbb{R}} |yz|^{H-\frac{3}{2}} f(y+z+x)dy dz,
\end{align*}
this completes the proof.
\fproof

Results from previous lemmas can be gathered in the extended equality which is given in the following lemma:

\begin{lemma} \label{1.3} Let  $\varphi,\psi:\mathbb{R}\rightarrow \mathbb{R}^{+}$. Then
	\begin{align} \label{eq l21}
	\int_{\mathbb{R}}\varphi(t)	M^2\psi(t)dt= \int_{\mathbb{R}}\psi(t)	M^2\varphi(t)dt=\int_{\mathbb{R}} M\varphi(t) M\psi(t) dt.
	\end{align}
\end{lemma}

\dproof 
From Lemma \ref{l1} it follows that
$$\int_{\mathbb{R}}\varphi(t)	M^2\psi(t)dt=\int_{\mathbb{R}}\varphi(t)C_H\int_{\mathbb{R}}\int_{\mathbb{R}}|yz|^{H-\frac{3}{2}}\psi(y+z+t)dydzdt.$$
If we substitute $ y+z+t=s,$ we obtain
\begin{eqnarray}
&&\int_{\mathbb{R}}\varphi(t)	M^2\psi(t)dt=C_H\int_{\mathbb{R}}\int_{\mathbb{R}}\int_{\mathbb{R}}|yz|^{H-\frac{3}{2}}\varphi(s-y-z)\psi(s)dydz{ds}\nonumber\\
&&\phantom{\int_{\mathbb{R}}\varphi(t)	M^2\psi(t)dt}=\int_{\mathbb{R}}C_H\int_{\mathbb{R}}\psi(s)\int_{\mathbb{R}}|yz|^{H-\frac{3}{2}}\varphi(s-y-z)dydzds\nonumber\\
&&\phantom{\int_{\mathbb{R}}\varphi(t)	M^2\psi(t)dt}=\int_{\mathbb{R}}\psi(s)	M^2\varphi(s)ds.
\end{eqnarray}
This proves the first identity. A similar argument proves the second identity.
\fproof

\subsection{The fractional Brownian motion}

We can now give our definition of fractional Brownian motion:

\begin{definition}{(Fractional Brownian motion)}\\
For $t \in \mathbb{R}$ define
\begin{equation} \label{eq4.5a}
B_t^{(H)}:= B^{(H)}(t):= B^{(H)}(t,\omega):=\langle \omega, M_t(\cdot) \rangle ; \quad \omega \in \Omega 
\end{equation}
where we for simplicity of notation put 
\begin{align}
M_t (x):=M(\chi_{[0,t]})(x)=(C_H)^{-1} \Big(\frac{t-x}{|t-x|^{\frac{3}{2}-H}}+\frac{x}{|x|^{\frac{3}{2}-H}}\Big);\quad t,x \in \mathbb{R}.
\end{align}

See (A7), p. 324 in \cite{EvdH}.
Then $B^{(H)}$ is a Gaussian process, and  we have 
\begin{equation}
\E[B^{(H)} (t)]=B^{(H)}(0)=0,
\end{equation}
and 
\begin{equation}\label{CH}
\E[B^{(H)}(s) B^{(H)}(t)]= C_H[ |t|^{2H}+ |s|^{2H}-|t-s|^{2H}]; \quad s,t \in \mathbb{R}.
\end{equation}
(See \cite{EvdH},(A10).)

Hence 
\begin{align}
\E[B^{(H)}(t)^2]= 2 C_H t^{2H}=:\sigma^2,
\end{align}\label{2.19}
and from this it follows,, since $B^{(H)}(t)$ is Gaussian,
\begin{align}
\E[(B^{(H)}(t))^4]= 3 \sigma^4 = 12 C_H^2 t^{4H}.
\end{align} \label{moment}


By the Kolmogorov continuity theorem the process $B^{(H)}$ has a continuous version, which we will also denote by $B^{(H)}.$  This process satisfies all the conditions of being an fBm given in Definition \ref{fBm} and we will use this as our
 \emph{ fractional Brownian motion (fBm) with Hurst parameter $H$}.
\end{definition}

\begin{remark}{\bf{(A remark about filtrations)}}
We let $\mathcal{F}=\{\mathcal{F}_t\}_{t\geq 0}$ denote the filtration generated by Brownian motion $B(t,\omega)$. This means that, for all $t$,  $\mathcal{F}_t$ is the sigma-algebra generated by the random variables $B(s,\cdot); s \leq t$. Similarly $\mathcal{F}^{(H)}=\{\mathcal{F}_t^{(H)}\}_{t\geq 0}$ denotes the filtration generated by $B^{(H)}(t,\omega)$.
Recalling that $B(t,\omega)=\langle \omega, \chi_{[0,t]} \rangle, $ while $B^{(H)}(t,\omega)=\langle \omega,M_t \rangle$, we see that the two filtrations are not identical. In particlar, $B^{(H)}(t)$ is not $\mathcal{F}$-adapted.\\
This is in contrast with the situation for some other constructions of fBm, where the filtration of $B(t)$ and $B^{(H)}(t)$ are the same. See e.g. \cite{BHOZ}.

\end{remark}
\subsection{Fractional white noise calculus}
Using that
\begin{equation}
B^{(H)}(t)=\int_{\mathbb{R}} \chi_{[0,t]}(s) dB^{(H)}(s)= \int_{\mathbb{R}}M_t(s)dB(s)=\langle \omega, M_t \rangle.
\end{equation}
we see that if $f$ is a simple function of the form $f=\sum c_k \chi_{[t_k, t_{k+1}]}$ then
\begin{align}
\int f(t) dB^{(H)}(t)&= \sum c_k (B^{(H)}(t_{k+1} )- B^{(H)}(t_k))=\sum c_k \langle \omega, M_{t_{k+1}} - M_{t_k} \rangle \nonumber\\
&=\langle \omega, \sum c_k (M_{t_{k+1} }-M_{t_k} \rangle 
=\langle \omega,M(\sum c_k \chi_{[t_k,t_{k+1}]}\rangle \nonumber\\
&= \langle \omega,Mf \rangle,
\end{align}
where $\chi$ denotes the indicator function.
From this it follows by approximating $f \in \mathcal{S}$ with step functions that
\begin{equation}
 \int_{\mathbb{R}} f(t) dB^{H}(t)=\langle \omega,Mf \rangle =\int_{\mathbb{R}} Mf(t) dB(t)
\end{equation}
for all $f \in \mathcal{S}$.

We know that the action of $\omega \in \mathcal{S}'(\mathbb{R})$ extends from $\mathcal{S}(\mathbb{R})$ to $L^2(\mathbb{R}),$ and that if $f \in \mathcal{S}$ then $Mf \in L^2(\mathbb{R}).$ Therefore we can define $M$ on $\mathcal{S}'(\mathbb{R})$ by setting
\begin{equation}
\langle M\omega ,\varphi \rangle= \langle \omega, M\varphi \rangle; \varphi \in \mathcal{S}(\mathbb{R}), \omega \in \mathcal{S}'(\mathbb{R}).
\end{equation}
Using this, we can write
\begin{align}
B^{(H)}(t,\omega)=\langle \omega, M_t \rangle= \langle M\omega, \chi_{[0.t]} \rangle =B(t,M\omega).
\end{align}

We define the \emph{fractional white noise} 
$W^{H}(t)=W^{H}(t,\omega)$ by
\begin{align}
W^{H}(t,\omega)=\frac{d B^{(H)}(t,\omega)}{dt} \text{ (derivative in } (\mathcal{S})^{*} ).
\end{align}
By the above it follows that 

\begin{align}
W^{H}(t,\omega)=\frac{d B(t,M\omega)}{dt}=W(t,M\omega),
\end{align}
 where 
\begin{align}
W(t)=W(t,\omega)=\frac{d B(t,\omega)}{dt} \text{ is the classical white noise in } (\mathcal{S})^{*}.
\end{align}
For $0 \leq T \leq \infty$ let $\mathcal{L}_H^2[0,T]$ denote the set of (measurable) stochastic processes $\varphi(t)=\varphi(t,\omega); \ t \geq 0, \omega \in \Omega$ such that 
\begin{align}
\E[\int_0^T (M \varphi)^2(t)dt] < \infty.
\end{align}
Then we define the \emph{Wick-Ito-Skorohod} (WIS)-integral  of $\varphi \in \mathcal{L}_H^2[0,T]$  as follows:

\begin{definition}  \label{wis}
(The  Wick--It\^o--Skorohod (WIS) integral) \\
Let  $\varphi : \mathbb{R} \rightarrow (\mathcal{S})^{*}$ be such that $\varphi (s)  \diamond W^{H}(s) $ is $ds$-integrable in $ (\mathcal{S})^{*}$. Then we say that  $\varphi$ is
$dB^{(H)}$-integrable and we define the Wick--It\^o--Skorohod (WIS) integral of $\varphi (\cdot) = \varphi (\cdot, \omega)$ with respect to $B^{(H)}$ by
\begin{align}
\int_{\mathbb{R}} \varphi(s) dB^{H}(s)= \int_{\mathbb{R}} \varphi(s) \diamond W^{H}(s) ds.
\end{align} 
\end{definition}

We recall the following theorem from \cite{EvdH}:
\begin{theorem}
Let $\varphi \in \mathcal{L}_H^2[0,T].$ \\
Then
\begin{equation} 
\int_{\mathbb{R}} \varphi(s) dB^{(H)}(s)=\int_{\mathbb{R}} M\varphi(s) dB(s)=\int_{\mathbb{R}} M\varphi(s) \diamond W(s) ds,
\end{equation} 
where $\diamond$ denotes the Wick product, the second integral is the classical Wick-Ito-Skorohod (WIS)-integral for $B$ and the last integral is the Lebesgue integral in the Hida space $(\mathcal{S})^{*}$ of stochastic distributions.
\end{theorem}

The next lemma states that the expectation of a WIS-integral of fractional Brownian motion is zero:

\begin{lemma}  Let $\varphi \in \mathcal{L}_H^2[0,T]$. Then
	$$\E\left[\int_0^{T}  \varphi(s) dB^{(H)}(s) \right]=0.$$
\end{lemma}
\dproof This follows from Definition \ref{wis} and the fact that  the expectation of a Wick product is the product of the expectations. Note that this holds also when the two factors are stochastic and not necessarily independent.
\fproof

Although this is not used later, we mention the following relation between the Wick product $\diamond$ and the ordinary product:
\begin{lemma}\cite{EvdH}
assume that  $ f, g \in L_H^2[0,T]$ (and hence $f,g$ are \emph{deterministic}).Then
\begin{align} \label{iso}
&(\int_0^t f(s) dB^{(H)}(s)) \diamond (\int_0^t g(s) dB^{(H)}(s))\nonumber\\
&=(\int_0^t f(s) dB^{(H)}(s)) (\int_0^t g(s) dB^{(H)}(s)) -(f,g)_{H}.
\end{align}
\end{lemma}

We need the following It\^o formula  for fractional WIS processes: 
\begin{theorem}\cite{BO} \label{tWIS}{\bf{(An It\^o formula for fractional WIS processes)}}\newline
Let $H > \frac{1}{2}$. Let $\alpha$ be a measurable process such that $\E[\int_0^t |\alpha(s)|ds] <\infty$ for all $t \geq 0$ and let $\beta$ be an $\mathcal{F}-$adapted, cadlag WIS-integrable process. Suppose $X(t)$ is a \emph{fractional WIS process} of the form
\begin{align}
dX(t)=\alpha(t)dt + \beta(t) dB^{(H)}(t);\quad X(0)=x \in \mathbb{R},
\end{align}
Let $f:[0,T] \times\mathbb{R} \mapsto \mathbb{R}$ be in  $C_{1,2}$ and put $Y(t)=f(t,X(t))$. Then
\begin{align}
dY(t)&=\frac{\partial}{\partial t}f(t,X(t)) dt +\frac{\partial}{\partial x}f(t,X(t)) dX(t) \nonumber\\
&+ \frac{ \partial^2}{\partial x^2} f(t,X(t))\beta(t) M^2(\beta \chi_{[0,t]})(t)dt.
\end{align}
\end{theorem}
\dproof
Since $\beta$ is adapted we have $MD_t^{(H)} \beta(s)=D_t \beta(s)=0$ for all $s < t$. Therefore this result follows as a special case of Theorem 2.11 in \cite{BO} .
\fproof



An important special case of the (2-dimensional) It\^o formula is the following  product rule:
\begin{lemma} \label{add2} (Product rule)
Let $\alpha_i,\beta_i; i=1,2$ be as in Theorem \ref{tWIS} and let
	
	$$dX_i=a_i(t)dt+b_i(t)dB^{(H)}(t), \quad i=1,2$$
	with $X_i(0)=x_i.$ Then the following holds:
	\begin{eqnarray}
	&&\E[X_1(t)X_2(t)]=x_1x_2+\E\Big[\int_0^t\Big( X_1(s)a_2(s)+X_2(s)a_1(s)\Big)ds\Big]\label{8}\\
	&&\phantom{E[X_1(t)X_2(t)]=}+\int_0^t\Big(b_1(s)M^2(b_2\chi_{[0,t]})(s)+b_2(s)M^2(b_1\chi_{[0,t]})(s)\Big) ds \Big]. \nonumber
	\end{eqnarray}
\end{lemma}
	
\section{An $L^2$-estimate for WIS-integrals}	
	The following theorem gives a fundamental $L^2-$  estimate for WIS-integrals with respect to fBm. We emphasise that the estimate is not an $L^2$ estimate of the usual type, because we are integrating only up to $t<\infty$ and the $L^2$ norm of the integral depends on \emph{both} the $L^2$ norm of the integrand and $t$. Thus the estimate blows up if we let $t$ go to infinity. Nevertheless, as we shall see, this estimate is sufficient for our study of SDEs up to a finite time $t$  in the next section.
	
	\begin{theorem} \label{inequalityt}
	Let  $\tfrac{1}{2}<H<1$ and let $\varphi(t)=\varphi(t,\omega)$   be an $\mathcal{F}$-adapted  stochastic process such that
	
	$$\E\left[\int_0^T \varphi^2(s)dt\right]<\infty.$$
	
	Define
	
	\begin{equation}
	X(t)=\int_0^t  \varphi(s) dB^{(H)}(s), \qquad 0\leq t \leq T,
	\end{equation}
	where the integral is interpreted in the WIS sense.
	
	Then
	\begin{equation}
	\E[X^2(t)]\leq K(H) \E\left[\int_0^t \varphi^2(s)dt\right] t^{2H-1}, \qquad t\leq T. \label{inequality} 
	\end{equation}
where
\begin{align}
K(H)=4C_H^2 \left( \frac{\pi 5^{2H-1}}{2H-1}+ \frac{2}{1-H}\right). \label{K}
\end{align}	
	\end{theorem}
	

\begin{remark}
It is an interesting question what is the best constant in (\ref{inequality}).
\end{remark}

\begin{remark}
Sufficient conditions (in terms of Malliavin derivatives) for the existence in $L^2(\P)$ of the WIS integral with respect to $B^{(H)}$ can be found in Bender \cite{Be1}, Corollary 3.5.\ .
For general $(\mathcal{S})^*$-valued integrands, conditions for the existence of the WIS-integral as an element of $(\mathcal{S})^*$ are given by Mishura \cite{M}, Theorem 2.3.1, using an inequality of Vaage \cite{V}.
\end{remark}

\noindent Proof of Theorem \ref{inequalityt}:\\
By splitting $\varphi$ into its positive and negative part  and multiplying the estimate by a factor of 2, we see that we may assume that $\varphi \geq 0$.
 By the product rule (Lemma \ref{add2} ) we have
$$A(t):=\E[X^2(t)]= \E\left[\int_0^t 2 \varphi(s)M^2(\varphi \chi_{[0,t]})(s) ds\right].$$
Also, by Lemma \ref{l1}
$$M^2f(x)=C_H^2\int_{\mathbb R} \int_{\mathbb R}|yz|^{H-\frac{3}{2}}f(x+y+z)dy dz.$$
Hence, using the Cauchy-Schwarz inequality twice we obtain

\begin{eqnarray*}
&&A(t)=2C_H^2\E\left[\int_0^t \left( \int_{\mathbb R} \int_{\mathbb R}|yz|^{H-\frac{3}{2}} \varphi(s) \varphi(s+y+z)  \chi_{[0,t]}(s+y+z) dy dz\right)ds\right]\\
&&\phantom{A(t)}=2C_H^2\left[ \int_{\mathbb R} \int_{\mathbb R}\left( \int_0^t \E ( \varphi(s) \varphi(s+y+z)  \chi_{[0,t]}(s+y+z) ds\right) |yz|^{H-\frac{3}{2}} dy dz\right]\\
&&\phantom{A(t)}\leq  2C_H^2\Big[ \int_{\mathbb R} \int_{\mathbb R}\left( \int_0^t \E [ \varphi^2(s) ]\chi_{[0,t]}(s+y+z)ds\right)^{\frac{1}{2}} \\
&&\hskip 3cm  \times \Big( \int_0^t \E[\varphi^2(s+y+z) ] \chi_{[0,t]}(s+y+z) ds\Big) ^{\frac{1}{2}}|yz|^{H-\frac{3}{2}} dy dz\Big]\\
&&\phantom{A(t)}\leq 2 C_H^2\left\{\left[ \int_{\mathbb R} \int_{\mathbb R}\int_0^t \E [ \varphi^2(s) ]|yz|^{H-\frac{3}{2}}    \chi_{[0,t]}(s+y+z) dsdy dz\right] ^{\frac{1}{2}}\right.\\
&&\phantom{A(t)\leq}\times  \left.\left[ \int_{\mathbb R} \int_{\mathbb R}\left( \int_0^t \E[\varphi^2(s+y+z) ] \chi_{[0,t]}(s+y+z) ds\right)|yz|^{H-\frac{3}{2}} dy dz\right] ^{\frac{1}{2}}\right\}\\
&&\phantom{A(t)}=2 C_H^2 J_1^{\frac{1}{2}}J_2^{\frac{1}{2}},
\end{eqnarray*}
where
\begin{eqnarray*}
&&J_1:=\int_{\mathbb R} \int_{\mathbb R}\int_0^t \E[\varphi^2(s)]  \chi_{[0,t]}(s+y+z) ds|yz|^{H-\frac{3}{2}} dy dz,\\
&&J_2:=\int_{\mathbb R} \int_{\mathbb R} \int_0^t \E[\varphi^2(s+y+z)]  \chi_{[0,t]}(s+y+z) ds|yz|^{H-\frac{3}{2}} dy dz.
\end{eqnarray*}

In order to estimate the $dsdydz$-integral $J_1$, we split it up into a $ds$-integral and a $dydz$-integral, as follows:\\
Note that
$$0\leq s \leq t \quad \& \quad s+y+z\in[0,t] \Rightarrow y+z \in[-t,t],$$

It follows that
\begin{align}
J_1\leq \int_0^t \E[\varphi^2(s)] ds\int_{\mathbb R} \int_{\mathbb R} \chi_{[-t,t]}(y+z) |yz|^{H-\frac{3}{2}} dy dz.
\end{align}
This $dydz$-integral has a singularity at $(0,0)$ and at infinity. To handle these singularities we divide the area of integration into 3 parts, namely a neighbourhood of $(0,0)$ and two neighbourhoods of infinity, as follows:
\begin{align}
\int_{\mathbb R} \int_{\mathbb R} \chi_{[-t,t]}(y+z) |yz|^{H-\frac{3}{2}} dy dz
\leq (a_1+a_2+a_3),
\end{align}
where

\begin{eqnarray*}
&&a_1=\int_{\mathbb R} \int_{\mathbb R} |yz|^{H-\frac{3}{2}} \chi_{[-3t,3t]\times[-4t,4t]}(y,z) dy dz,\\
&&a_2=\int_{\mathbb R} \int_{\mathbb R}  |yz|^{H-\frac{3}{2}} \chi_{(-\infty,-3t]}(y)\chi_{[-t,t]}(y+z)  dy dz\\
&&a_3=\int_{\mathbb R} \int_{\mathbb R}  |yz|^{H-\frac{3}{2}} \chi_{(3t,\infty)}(y) \chi_{[-t,t]}(y+z)  dy dz.
\end{eqnarray*}

Note that
$$ (y,z) \in [-3t,3t] \times [-4t,4t] \Rightarrow y^2 + z^2 \leq 25 t^2 = (5t)^2.$$

Hence, by introducing polar coordinates, we can estimate the integral $a_1$ as follows:

\begin{align*}
a_1\leq \int_{0}^{2\pi} \int_0^{5 t} |r^2\sin \alpha\cos\alpha|^{H-\frac{3}{2}} r dr d\alpha
\leq \frac{2\pi 5^{2H-1} }{2H-1} t^{2H-1}.
\end{align*}

Next, note that
$$y\leq -3t \quad \& \quad  y+z\geq -t \quad  \Rightarrow z  \geq 2t.$$
Therefore,
\begin{eqnarray*}
&&a_2=\int_{\mathbb R} \int_{\mathbb R}  |yz|^{H-\frac{3}{2}} \chi_{(-\infty,-3t]}(y)\chi_{[-t,t]}(y+z)  dy dz\\
&&\phantom{a_2}\leq\int_{2t}^{+\infty}\left(\int_{-t-z}^{t-z}|y|^{H-\frac{3}{2}}dy\right)|z|^{H-\frac{3}{2}}dz\\
&&\phantom{a_2}=\int_{2t}^{+\infty}\left(\frac{-1}{H-\frac{1}{2}}(-y)^{H-\frac{1}{2}}\Big|_{-t-z}^{t-z}\right)|z|^{H-\frac{3}{2}}dz\\
&&\phantom{a_2}=\frac{1}{H-\frac{1}{2}}\int_{2t}^{+\infty}\left((z+t)^{H-\frac{1}{2}}-(z-t)^{H-\frac{1}{2}}\right)|z|^{H-\frac{3}{2}}dz\\
&&\phantom{a_2}\leq\frac{1}{H-\frac{1}{2}}\int_{2t}^{+\infty}\left(H-\frac{1}{2}\right)(z-t)^{H-\frac{3}{2}}|z|^{H-\frac{3}{2}}2tdz\\
&&\phantom{a_2}\leq2t\int_{2t}^{+\infty} (z-t)^{2H-3}dz\\
&&\phantom{a_2}=2t\frac{(z-t)^{2H-2}}{2H-2}\Big|_{2t}^{+\infty}\\
&&\phantom{a_2}=\tfrac{1 }{1-H}t^{2H-1}.
\end{eqnarray*}
	
	Similarly we obtain
	
	$$a_3\leq \tfrac{1 }{1-H}t^{2H-1}.$$
	
	In order to estimate $J_2$ we make a change of variables:
	
	\begin{equation*} 
	\left\{ 
        \begin{array}{l} 
            u=s+y+z\\
            v=y\\
            w=z\\
        \end{array}
        \right. 	\Leftrightarrow  \begin{array}{l} 
            s=u-v-w=: X(u,v,w)\\
            y=v:=Y(u,v,w)\\
            z=w=:Z(u,v,w).\\
        \end{array}
    \end{equation*}

    Then the Jacobian is \\
    \begin{equation}
    \begin{vmatrix}
\frac{\partial X}{\partial u} & \frac{\partial X}{\partial v} & \frac{\partial X}{\partial w} \\ 
\frac{\partial Y}{\partial u} & \frac{\partial Y}{\partial v} & \frac{\partial Y}{\partial w} \\ 
\frac{\partial Z}{\partial u} & \frac{\partial Z}{\partial v} & \frac{\partial Z}{\partial w} 
\end{vmatrix} =   
    \begin{vmatrix}
1 & -1 & -1 \\ 
0 & 1 & 0\\ 
0 & 0 & 1
\end{vmatrix}=1.
\end{equation}	
\vskip 1cm 
	So we get
	
	\begin{equation*}
	    J_2=\int_{\mathbb R} \int_{\mathbb R}\int_0^t \E[\varphi^2(u)]  \chi_{[0,t]}(u) ds|vw|^{H-\frac{3}{2}} dv dw,
	\end{equation*}
	which is of the same type as $J_1$.
	
	We conclude that
	
	$$ \E[X^2(t)]\leq 2 C_H^2 \left( \frac{\pi 5^{2H-1}}{2H-1}+ \frac{2}{1-H}\right)E\left[\int_0^T \varphi^2(s)dt\right] t^{2H-1}, \qquad t\leq T,$$
	which completes the proof.
\fproof

\section{Stochastic differential equations driven by fractional Brownian motion}

In this section we prove our main result, namely the existence and uniqueness of the solution of the following WIS-type SDE driven by fBm: 

\begin{theorem}\label{sde3}
Let $T>0$ be fixed and\\
$\alpha(\cdot,\cdot,\cdot)\colon[0,T]\times\RR\times \Omega \to\RR, \\
\beta(\cdot,\cdot,\cdot)\colon[0,T]\times\RR \times \Omega \to\RR$ and \\
$\sigma(\cdot,\cdot,\cdot)\colon[0,T]\times \R \times \Omega \to
\RR$\\
 be stochastic processes satisfying
\begin{equation}
|\alpha(t,x,\omega)|+|\beta(t,x,\omega)| + |\sigma(t,x,\omega)|\leq C(1+|x|)\;;\qquad x \in \RR,\;\label{assam.1}
t\!\in\![0,T]
\end{equation}
and
\begin{align}
&|\alpha(t,x,\omega)-\alpha(t,y, \omega)|+|\beta(t,x,\omega)-\beta(t,y, \omega)| +|\sigma(t,x,\omega)-\sigma (t,y,\omega)|\nonumber\\
&\leq D|x-y|\,;\quad
x,y\!\in\!\RR,\, t\!\in\![0,T] \;\label{assam.2}
\end{align}
for some constants $C,D$ not depending on $t,x,y$.

 Let
$Z$ be a random variable which is independent of both $\mathcal{F}$ and $\mathcal{F}^{(H)}$
and such that
$$
E[|Z|^2]<\infty\;.
$$
Then the stochastic differential equation
\begin{align}  \label{SDE}
dX_t&=\alpha(t,X_t, \omega)dt+\beta(t,Y_t,\omega) dB(t)+\sigma(t,Y_t,\omega)dB^{(H)}_t\\
&Y_t=\E[X_t | \mathcal{F}_t];,\qquad 0\!\leq\!t\leq\!T,\,
X_0=Z
\end{align}
has a unique solution $X_t(\omega)$ with the property that
\begin{equation}
\E\bigg[ \int\limits_0^T |X_t|^2dt\bigg]<\infty\;.
\end{equation}
\end{theorem}

\dproof
We proceed as in the proof of Theorem 5.2.1 in \cite{O}:\\
We may assume that $\beta=0$, since the handling of the Brownian motion integral is well-known.
\medskip

{\it Uniqueness.} \quad Let
$X_1(t,\omega)=X_t(\omega)$ and $X_2(t,\omega)=\widehat{X}_t(\omega)$ be two
solutions of \eqref{SDE} with initial values $Z,\widehat{Z}$ respectively, i.e.
$Z(\omega):=X_1(0,\omega), \widehat{Z}(\omega):=X_2(0,\omega), \omega\in\Omega$.

Put $\rho(s)=\alpha(s,X_s,\omega)-\alpha(s, \widehat{X}_s,\omega)$ and $\gamma(s)=\sigma(s,Y_s)-\sigma(s,\widehat{Y}_s)$, where $\widehat{Y}_s=\E[\widehat{X}_s | \mathcal{F}_s]$. Then by Theorem \ref{inequalityt}  we get (in the following we suppress the argument with respect to time and $\omega$),
     \begin{eqnarray*}
&&\hskip -0.6 cm\E[|X_t-\widehat{X}_t|^2]=\E\bigg[\bigg(
    Z-\widehat{Z}+\int\limits_0^t \rho\,ds + \int_0^t \gamma(s) dB^{(H)}(s)\Big)^2 \bigg] \\
&&\hskip -0.6 cm\phantom{\E[|X_t-\widehat{X}_t|^2]}\leq 3\E[|Z-\widehat{Z}|^2]+3\E\bigg[\bigg(\int\limits_0^t \rho\,ds \bigg)^2\bigg]+3\E\Big[\Big(\int_0^t \gamma(s) dB^{(H)}(s)\Big)^2 \Big]\\
&&\hskip -0.6 cm\phantom{\E[|X_t-\widehat{X}_t|^2]}\leq 3\E[|Z-\widehat{Z}|^2]+3t\E\Big[\int\limits_0^t \rho^2
     ds\Big] + 3K(H) t^{2H-1} \E\Big[\int_0^t \gamma(s)^2 ds\Big],
     \end{eqnarray*}
and further
\begin{eqnarray*}
&&\E[|X_t-\widehat{X}_t|^2]\leq 3\E[|Z-\widehat{Z}|^2]+3tD^2\int\limits_0^t
     \E[|X_s-\widehat{X}_s|^2]ds\\
     &&\phantom{\E[|X_t-\widehat{X}_t|^2]\leq}+3K(H) t^{2H-1}D^2\int\limits_0^t
     \E[|Y_s-\widehat{Y}_s|^2]ds)\\
       &&\phantom{\E[|X_t-\widehat{X}_t|^2]}\leq3\E[|Z-\widehat{Z}|^2]+3tD^2\int\limits_0^t
     \E[|X_s-\widehat{X}_s|^2]ds\\
  &&\phantom{\E[|X_t-\widehat{X}_t|^2]\leq}+3K(H) t^{2H-1}D^2\int\limits_0^t
     \E[|X_s-\widehat{X}_s|^2]ds),
\end{eqnarray*}
\noindent where $D$ is from (\ref{assam.2}).

So the function
$$
w(t)=\E[|X_t-\widehat{X}_t|^2]\;;\qquad 0\leq t\leq T
$$
satisfies
\begin{eqnarray}
&& w(t)\leq F+A\int\limits_0^t w(s)ds\;,   \\
&&\hskip-4.3cm \mbox{where $\,F=3\E[|Z-\widehat{Z}|^2]$ and 
       $\,A=3D^2(T+K(H)T^{2H-1})$}\;. \nonumber
\end{eqnarray}
By the Gronwall-Bellman inequality  (see \cite{bainov}) we conclude that
\begin{equation}
w(t)\leq F\exp(At)\;.
\end{equation}
Now assume that $Z=\widehat{Z}$. Then $F=0$ and so $w(t)=0$ for all $t\geq 0$. 
Hence
$$
X_t= \widehat{X}_t \text{ a.s. for all } t \in [0,T],
$$
and the uniqueness is proved.\\

{\it Existence.} \quad The proof of the existence is similar to the familiar existence proof for
ordinary differential equations. For completeness we give the details: \\

Without loss of generality we may assume that $T\leq 1$.

\smallskip

Define $X_t^{(0)}=X_0$ and $X_t^{(k)}=X_t^{(k)}(\omega)$ inductively as 
follows
\begin{equation}
X_t^{(k+1)}=X_0+\int\limits_0^t
\alpha(s,\omega, X_s^{(k)})ds+\int\limits_0^t
\sigma(s,Y_s^{(k)}, \omega)dB_s^{(H)}\;.
\end{equation}
Then, similar computation as for the uniqueness above gives
$$
\E[|X_t^{(k+1)}-X_t^{(k)}|^2]\leq 3D^2(T+K(H)T^{2H-1})\int\limits_0^t
\E[|X_s^{(k)}-X_s^{(k-1)}|^2]ds\;,
$$
for $k\geq 1$, $t\leq T$
and
$$
\E[|X_t^{(1)}-X_t^{(0)}|^2]\leq 2C^2t^2
         \E[(1+|X_0|)^2] \leq A_1t,
$$
where the constant $A_1$ only depends on $C,T$ and $E[|X_0|^2]$. So by
induction on $k$ we obtain
\begin{equation}
\E[|X_t^{(k+1)}-X_t^{(k)}|^2]\leq\frac{A_2^{k+1}t^{k+1}}{(k+1)!}\;; 
\qquad k\geq 0,\;t\in[0,T]
\end{equation}
for some suitable constant $A_2$ depending only on $C,D,T$ and $\E[|X_0|^2]$. 

Hence, if $\lambda$ denotes Lebesgue measure on $[0,T]$ and $m>n\geq 0$ we get
\begin{eqnarray}
&& \big\| X_t^{(m)}-X_t^{(n)}\big\|_{L^2(\lambda\times \P)}=
  \Big\| \sum_{k=n}^{m-1} X_t^{(k+1)}-X_t^{(k)} 
        \Big\|_{L^2(\lambda\times \P)} \nonumber \\
&&\quad \leq \sum_{k=n}^{m-1} \big\| X_t^{(k+1)}
      -X_t^{(k)}\big\|_{L^2(\lambda\times \P)}
   =\sum_{k=n}^{m-1} \Big( \E \Big[ \int\limits_0^T
     |X_t^{(k+1)}-X_t^{(k)}|^2 dt \Big] \Big)^{1/2} \nonumber \\
&&\quad \leq \sum_{k=n}^{m-1} \Big( \int\limits_0^T
   \frac{A_2^{k+1}t^{k+1}}{(k+1)!}dt\Big)^{1/2}
   =\sum_{k=n}^{m-1} \Big(
   \frac{A_2^{k+1}T^{k+2}}{(k+2)!}\Big)^{1/2}\to 0,
\end{eqnarray}
as $m,n\to\infty$.

Therefore $\{X_t^{(n)}\}_{n=0}^\infty$ is a Cauchy sequence in
$L^2(\lambda\times \P)$. Hence $\big\{X_t^{(n)}\big\}_{n=0}^\infty$ is
convergent in
$L^2(\lambda\times \P)$. Define
$$
X_t:=\lim_{n\to\infty} X_t^{(n)}\qquad \hbox{(limit in $L^2(\lambda\times
\P)$)}.
$$
It remains to prove that $X_t$ satisfies \eqref{SDE}:

For all $n$ and all $t\in[0,T]$ we have
$$
X_t^{(n+1)}=X_0+\int\limits_0^t \alpha(s,X_s^{(n)},\omega)ds
   +\int\limits_0^t \sigma(s, Y_s^{(n)},\omega)dB_s^{(H)}\;.
$$
Now let $n\to\infty$. Then by the estimates above we get that
$$
\int\limits_0^t \alpha(s, X_s^{(n)},\omega)ds\to
\int\limits_0^t \alpha(s, X_s, \omega)ds\qquad \hbox{in $\;L^2(\P)$}.
$$
and
$$
\int\limits_0^t \sigma(s, Y_s^{(n)},\omega)dB^{(H)}(s)\to
\int\limits_0^t \sigma(s, Y_s, \omega)dB^{(H)}(s)\qquad \hbox{in $\;L^2(\P)$}.
$$
We conclude that for all
$t\in[0,T]$ we have
\begin{equation}
X_t=X_0+\int\limits_0^t \alpha(s, X_s,\omega)ds
  +\int\limits_0^t \sigma(s, Y_s, \omega)dB_s^{(H)}\quad \hbox{a.s.}
\end{equation}
This completes the proof.
\fproof

\begin{remark}
Note that even if we assume that $Z$ is constant and $\alpha(t,x,\omega)=\alpha(t,x)$ and $\sigma(t,\omega)=\sigma(t)$ do not depend on $\omega$, our proof does not give that the solution $X_t$ is $\mathcal{F}$- adapted. \\
Moreover, it is not clear under what conditions we would have that, for almost all t, the map $t \mapsto X_t$ is continuous a.s. A proof of this seems to require a Burkholder-Davis-Gundy type inequality for WIS-integrals w.r.t. fBm. This is a topic for future research.
\end{remark}


\end{document}